\documentclass[12pt, a4wide]{article}
\title{Hironaka's characteristic polygon and effective resolution of surfaces}
\author{Piedra, R.\footnote{Supported by Universidad de Sevilla and MCyT
BFM2000-1523} \and Tornero, J.M.\footnote{Supported by FQM 218 and
MTM2004--07203--C02--01.} }
\date{August, 2005}

\newcommand{\NN}{{\bf N}}

\newcommand{\VS}{\vspace{.5cm}}

\newcommand{\vs}{\vspace{.15cm}}

\newcommand{\df}{\noindent {\bf Definition.-- }}

\newcommand{\thr}{\noindent {\bf Theorem.-- }}

\newcommand{\obs}{\noindent {\bf Remark.-- }}

\newcommand{\cS}{{\cal S}}

\DeclareMathAlphabet{\frak}{U}{euf}{m}{n}
\SetMathAlphabet{\frak}{bold}{U}{euf}{b}{n}

\newcommand{\fp}{{\frak p}}

\begin{document}

\maketitle

\abstract{Hironaka's concept of characteristic polyhedron of a
singularity has been one of the most powerful and fruitful ideas
of the last decades in singularity theory. In fact, since then
combinatorics have become a major tool in many important results.
However, this seminal concept is still not enough to cope with
some effective problems: for instance, giving a bound on the
maximum number of blowing--ups to be performed on a surface before
its multiplicity decreases. This short note shows why such a
bounding is not possible, at least with the original definitions.}


\section{Introduction}

In this paper we will deal with embedded algebroid surfaces, that
is, schemes given by the spectrum of a ring $R = K[[X,Y,Z]]/(F)$,
where $K$ is an algebraically closed field and $F$ is a power
series of order $n>0$. Such an $F$ will be called an equation of
the surface and $n$ will be called the multiplicity of the
surface.

\vs

\thr ({\bf Levi--Zariski, algebroid version}) Let $\cS$ be an
algebroid embedded surface with normal crossing singularities.
Then if we sucessively blow up smooth equimultiple subvarieties of
maximal dimension, the multiplicity is dropped in a finite number
of steps.

\vs

The fact that a surface can be resolved by blowing up ``maximal''
centers was already investigated by Beppo Levi (\cite{L2})
and actually proved by Zariski (\cite{Zariski}) in characteristic
zero as a part of his proof of the resolution for three--dimensional
varieties. His results and techniques were paralleled by Abhyankar in
positive characteristic (\cite{Abh}).

\vs

It was Hironaka, however, who got a massive breakthrough by
considering for the first time the use of combinatorial tools in
singularity theory (\cite{H1,Bowdoin}), a technique who has made
possible to tackle a number of problems (see for instance
\cite{Cossart,WPG,HPG,HH} or, for a much more complete
information, the excellent survey \cite{17}).

\vs

As for surfaces is concerned, a combinatorial approach to the
Levi--Zariski theorem was pointed out by Hironaka in the
introduction of \cite{H1}, and partially developed in
\cite{Bowdoin}, although the arguments are not very clear in the
positive characteristic case. The original purpose of this paper
was sharing the combinatorial approach, using a somehow different
induction argument to produce an upper bound for the number of
blowing--ups that can be performed, following the Levi--Zariski
procedure, before the multiplicity drops. In the case of curves,
this was attached quite straightforwardly using the first
characteristic exponent (see, for instance, \cite{H1} and
\cite{ZZ}). However, we have found that already for surfaces this
bounding is not possible (at least using Hironaka's characteristic
polygon).

\vs

The authors are enormously grateful to Prof. J.L. Vicente, from
whom they learned the subject and who posed this question to the
last author some years ago. This paper is dedicated to him on his
1000000th (binary) birthday .

\section{Some technical set--up}

For the sake of completeness, we recall here well-known technical results
that will be of some help in the sequel.

\vs

Let $\cS$ be an embedded algebroid surface of multiplicity $n$,
$F$ an equation of $\cS$. After a change of variables, one can
take $F$ to the (Weierstrass) form
$$
F(X,Y,Z) = Z^n + \sum_{k=0}^{n-1} a_kZ^k, \mbox{ where } a_k(X,Y)  =
\sum_{i,j} a_{ijk} X^iY^j.
$$

To this situation a combinatorial object may be attached: note
$$
N_{\{X,Y,Z\}}(F) = \left\{ (i,j,k) \in \NN^3 \; | \; a_{ijk} \neq 0 \right\}
\cup \{ (0,0,n) \},
$$
where we will omit the subscript whenever the variables are clear
from the context. The Hironaka (or Newton, or
Newton--Hironaka,...) polygon of $F$ is
$$
\Delta_{\{X,Y,Z\}} (F) (\mbox{or } \Delta (F)) = \mbox{CH} \left( \bigcup_{a_{ijk}\neq 0} \left[  \left(
\frac{i}{n-k}, \frac{j}{n-k} \right) + \NN^2 \right]  \right),
$$
where CH stands for the convex hull. This object was already used
in the famous Bowdoin notes lectures by Hironaka (\cite{Bowdoin})
and it appeared in printed form for the first time in the
outstanding paper \cite{H1}, where it covered the surface case of
the much more general notion of characteristic polyhedron of a
singularity.

\vs

\obs If we allow $Z$ to vary, using changes of variable of the
type
$$
Z \; \longmapsto \; Z + \alpha(X, Y), \mbox{ with } \alpha \in K[[X,Y]] \mbox{ not a unit},
$$
we obtain a collection of polygons which has a minimal element in
the sense of inclusion (this is not obvious at all in positive
characteristic). This object was called by Hironaka the
characteristic polygon of the pair $(\cS, \{ X,Y\})$ (\cite{H1}), noted
$\Delta (\cS, \{X,Y\})$.

\vs

%
%
%
%
%

\vs

\obs For the case of characteristic zero or, more broadly, the case
where $n$ does not divide the characteristic of $K$, it is customary to
make the Tchirnhausen transformation,
$$
Z \; \longmapsto \; Z - \frac{1}{n} a_{n-1} (X,Y),
$$
which is a change of variables of the type considered by Hironaka.
Such an equation will be called a WT equation and it has many
interesting properties, which we will show. To begin with, in
these equations, a permitted (that is, equimultiple and smooth)
curve can be written in the form $\fp = (Z,G(X,Y))$.

\vs

\df A vertex $(P_1,P_2)$ of $\Delta (F)$ is called contractible if
there exists a change of variables $\varphi$
$$
Z \; \longmapsto \; Z + \alpha X^a Y^b, \mbox{ with }
\alpha \in K,
$$
such that
$$
\Delta (\varphi(F)) \subset \Delta (F) \setminus \{
(P_1,P_2) \}.
$$

Were this the case, $\varphi$ is called the contraction of the vertex
$(P_1,P_2)$.

\vs

\obs If we can apply the Tchirnhausen transformation, the
resulting equation has no contractible vertices. In fact, a vertex
$(a,b)$ is contractible if and only if it represents {\em all} the
monomials from $(Z + \alpha X^a Y^b)^n$ and this cannot happen
since $a_{n-1}(X,Y)=0$. As it will become obvious from the
equations associated to the different blowing--ups, this situation
will remain during the resolution process (at least, until a
multiplicity decrease happens). In classical terms, $Z=0$ is a
linear hypersurface with permanent maximal contact with the
surface $\cS$.

\vs

Hironaka proved in \cite{H1} (for arbitrary characteristic)
that all the vertices of $\Delta (F)$ are not contractible if and
only if $\Delta (F) = \Delta (\cS, \{X,Y\})$. From the previous
remark this is obvious in the case of WT equations.

\section{An interesting example}

We will try to prove that Hironaka's characteristic polygon does not
contain enough information in order to bound the resolution process,
even if we are interested only in a first multiplicity decreasing.

\vs

To see that, assume that $K$ has characteristic other than $3$ and
consider the surface $\cS$ defined by the equation
$$
F= Z^3 + X^mZ + (X-Y)^4, \mbox{ with } m\geq 19;
$$
\noindent which is, obviously, a WT equation and, hence, the
characteristic polygon of $\cS$ is given by

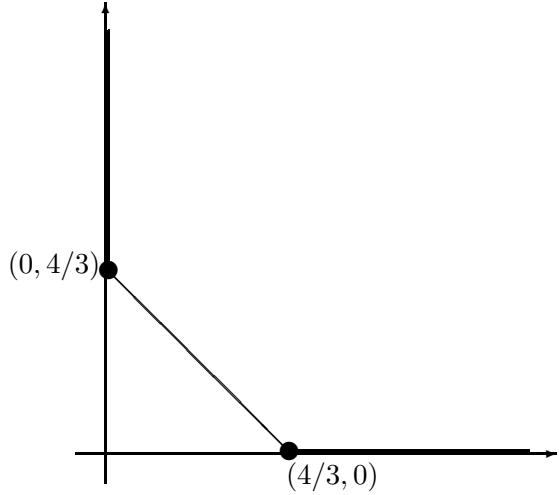
\vs

\unitlength=.4cm

\begin{picture}(20,20)(0,0)

\put(4,3){\vector(0,1){16}} \put(3,4){\vector(1,0){16}}

\put(4.1,10.1){\line(0,1){8}} \put(4.1,10.1){\line(1,-1){6}}
\put(10.1,4.1){\line(1,0){8}}

\put(4.075,10.075){\line(0,1){8}}
\put(4.075,10.075){\line(1,-1){6}}
\put(10.075,4.075){\line(1,0){8}}

\put(.75,10){$(0,4/3)$} \put(10,3){$(4/3,0)$}

\put(4.1,10.1){\circle*{.6}} \put(10.1,4.1){\circle*{.6}}

\end{picture}

\vs

We will first make, in absence of permitted curves, a quadratic
transformation (that is, blowing--up the origin) on the direction
$(1:1:0)$, giving as a result the surface $\cS_1$ defined by
$$
F_1 = Z^3 +X^{m-2}Z + XY^4.
$$

Next we will make two quadratic transformations on the direction
$(1:0:0)$ (no permitted curves in either surfaces) after which we
get $\cS_3$, defined by
$$
F_3 = Z^3 + X^{m-6}Z + X^3Y^4.
$$

We are hence forced now to perform a monoidal transformation
(blowing--up of a curve) centered on $(Z,X)$, as it is now
permitted. We get then $\cS_4$, defined by
$$
F_4 = Z^3 + X^{m-8}Z + Y^4.
$$

Now it is straightforward that, after three quadratic
transformations centered on $(1:0:0)$ and a monoidal
transformation center on $(Z,X)$, we should get a surface $\cS_8$
defined by
$$
F_8 = Z^3 + X^{m-16}Z + Y^4.
$$

Obviously this implies that it is not possible to get a bound for
the number of blowing--ups needed for decreasing the multiplicity
of $\cS$: as we change $m$ we get a family of surfaces with the
same Hironaka polygon but needing an arbitrarily large number of
blowing--ups to get a multiplicity loss.

\vs

\obs The key for this counterexample is the first quadratic
transformation; in fact it is easy to bound a resolution process
using uniquely monoidal transformations or quadratic
transformations centered in $(1:0:0)$ and/or $(0:1:0)$ (this is
connected with the Weak Hironaka's Polyhedra Game, solved by
Spivakovsky in \cite{WPG}).

However, as this example makes apparent, if we ever want to bound
the resolution process we need to be able to track much more
complicated relations between our parameters than the ones
considered by Hironaka. We hope this example sheds some light in
order to tackle this (for us) quite interesting problem.

\VS\VS

\noindent R. Piedra.

Departamento de \'Algebra.

Facultad de Matem\'aticas.

Universidad de Sevilla.

Apdo. 1160. 41080 Sevilla (Spain).

\noindent E--mail: \verb|piedra@algebra.us.es|

\vs

\noindent J.M. Tornero.

Departamento de \'Algebra.

Facultad de Matem\'aticas.

Universidad de Sevilla.

Apdo. 1160. 41080 Sevilla (Spain).

\noindent E--mail: \verb|tornero@algebra.us.es|

\end{document}